\documentclass[oneside,leqno,11pt]{article}
\usepackage{amssymb,amsmath,latexsym}

\def\gk{\mathfrak{k}}

\def\C{\mathbb{C}}

\def\E{\mathbb{E}}

\def\P{\mathbb{P}}

\def\R{\mathbb{R}}

\def\Z{\mathbb{Z}}

\def\cD{\mathcal{D}}

\def\cH{\mathcal{H}}
\def\cI{\mathcal{I}}

\def\cP{\mathcal{P}}

\newtheorem{theorem}[equation]{Theorem}

\newtheorem{lemma}[equation]{Lemma}

\newtheorem{proposition}[equation]{Proposition}

\newtheorem{remark}[equation]{Remark}

\title{Spherical Functions on Euclidean Space}

\date{20 September 2005}

\author{Joseph A. Wolf\,\footnote{
Research partially supported by NSF Grant DMS 04-00420.}}

\begin{document}
\maketitle

\begin{abstract}
We study special functions on euclidean spaces from the viewpoint of 
riemannian symmetric spaces.  Here the euclidean space $E^n = G/K$ where 
$G$ is the semidirect product $R^n \cdot K$ of the translation group with 
a closed subgroup  $K$ of the orthogonal group $O(n)$.  We give exact 
parameterizations of the space of $(G,K)$--spherical functions by a certain 
affine algebraic variety, and of the positive definite ones by a real form 
of that variety.  We give exact formulae for the spherical functions in 
the case where $K$ is transitive on the unit sphere in $E^n$

\end{abstract}

\section{Introduction} \label{sec1}
\setcounter{equation}{0}

Riemannian symmetric spaces carry a distinguished class of functions,
the {\em spherical function}, that generalize the notion of characters
on a commutative group.  In particular, if $X = G/K$ is a riemannian
symmetric space, then $L^2(X)$ is a continuous direct sum (direct
integral) of Hilbert
spaces $E_\varphi$ generated by positive definite spherical functions
$\varphi: X \to \C$, $G$ acts naturally on $E_\varphi$ by an
irreducible unitary representation $\pi_\varphi$, and the left regular 
representation of $G$ on $L^2(X)$ is a continuous direct sum of the 
$\pi_\varphi$.  This circle of ideas has been studied in detail for a
number of years.

With one ``small'' exception, the symmetric space $X$ determines a
standard choice of $(G,K)$, as follows: $G$ is the group $\cI(X)$ of all
isometries of $X$, or the subgroup generated by all geodesic symmetries
of $X$, or the identity component $\cI(X)^0$, and the latter is the group 
generated
by all products of an even number of geodesic symmetries; in any case $K$ 
is the isotropy subgroup of $G$ at a ``base point'' $x_0 = 1K$.  The
exception is when $X$ is of the form $X' \times \E^n$ where $\E^n$
is an euclidean space of dimension $n > 0$ and $X'$ has no euclidean factor.  
Then $G = G' \times J$ where $G'$ is determined by $X'$ as above and where
$J$ is any closed subgroup of the euclidean group $E(n)$ that contains the
translations.  There are infinitely many natural choices of $J$, and the 
spherical functions have only been studied in a serious way for two of those 
choices.  In this note we hope to close that gap by presenting a general
theory of spherical functions on euclidean spaces.  The main results are
Theorems \ref{maintheorem1} and \ref{ratparam} concerning spherical functions
in general and their parameterization by a certain affine variety,
Theorem \ref{ratparamreal} which parameterizes the positive definite
spherical functions, and Theorems \ref{catquottrans} and 
\ref{trans-case-sph-fns} concerning the case where the isotropy group
is transitive on the unit sphere in the tangent space.

\section{Background on Spherical Functions} \label{sec2}
\setcounter{equation}{0}

Spherical functions on a riemannian symmetric space $X = G/K$ can be 
characterized in a number of ways.  We recall the relevant characterizations
from \cite{Di}, \cite{H1}, \cite{H2} and \cite{W3}.

The basic definition, which does not use differentiability, is this.
A function $\varphi: G \to \C$ is spherical relative to $(G,K)$
if (i) it is continuous, (ii) it is bi--$K$--invariant (in other words,
is a $K$--invariant function on $X = G/K$), (iii) it is normalized
by $\varphi(1) = 1$, and (iv) if $f \in C_c(K\backslash G/K)$ there exists
$\lambda_f \in \C$ such that $f*\varphi = \lambda_f\,\varphi$.  Here
$C_c(K\backslash G/K)$ denotes the space of all compactly supported
bi--$K$--invariant continuous functions on $G$, in other words continuous
compactly supported $K$--invariant functions on $X$.

Spherical functions are also characterized as the continuous functions
$\varphi: G \to \C$, not identically zero, that satisfy the functional
equation: if $g_1, g_2 \in G$ then the integral
$\varphi(g_1)\varphi(g_2) = \int_K \varphi(g_1kg_2)\, d\mu_{_K}(k)$.
We will make serious use of the functional equation.

Finally, spherical functions $\varphi : G \to \C$ are characterized as 
the joint eigenfunctions 
of the algebra $\cD(G,K)$ of $G$--invariant differential operators on
$X$, normalized to take the value $1$ at the base point $x_0 = 1K$.
Note that the joint eigenvalue is an associative algebra homomorphism 
$\chi_\varphi : \cD(G,K) \to \C$.  We will need the fact that if
$\varphi$ and $\varphi'$ are spherical functions, and if the joint
eigenvalues $\chi_\varphi = \chi_{\varphi'}$, then $\varphi = \varphi'$.

The notion of induced spherical function mirrors the notion of induced
representation.  Let $Q \subset G$ be a closed subgroup such that
$K$ is transitive on $G/Q$, i.e. $G = KQ$, i.e. $G = QK$, i.e.
$Q$ is transitive on $G/K$.
Let $\zeta: Q \to \mathbb C$ be spherical for $(Q,Q\cap K)$.  The {\em
induced spherical function} is 
\begin{equation} \label{indrep}
[\text{Ind}_Q^G(\zeta)](g)
= \int_K \widetilde{\zeta}(gk) \, d\mu_{_K}(k)
\text{ where } \widetilde{\zeta}(kq) = \zeta(q)\Delta_{G/Q}(q)^{-1/2}.
\end{equation}
Here $\Delta_{G/Q}: Q \to \R$ is the quotient of modular functions,
$\Delta_{G/Q}(q) = \Delta_G(q)/\Delta_Q(q) = \Delta_Q(q)^{-1}$.

\section{Background on Euclidean Space} \label{sec3}
\setcounter{equation}{0}

Let $K$ be any closed subgroup of the orthogonal group $O(n)$ and consider
the semidirect product group $G = \R^n\cdot K$.  That gives an expression
of euclidean space as a riemannian symmetric coset space, $\E^n = G/K$.  
Here $K$ is identified with the isotropy subgroup of $G$ at the base point 
$(0,K)$ corresponding to a choice of origin.  The extreme cases are
$K = \{1\}$, where $G$ is just the group of euclidean translations of $\E^n$,
and $K = O(n)$ or $SO(n)$, where $G$ is the full euclidean group $E(n)$
or its identity component $E(n)^0$.  

The group $K$ acts on $\C^n$ by complexification of its natural action
on $\R^n$ as a subgroup of $O(n)$.  It preserves the $\C$--bilinear
form $b(\xi,\eta)$ that is the complex extension of the $O(n)$--invariant
inner product on $\R^n$.  If $\xi \in \C^n$ we have the quasicharacter
$\varphi_\xi : \R^n \to \C$ given by 
$\varphi_\xi(x) = e^{i b(x,\xi)}$.  In the identification of
$\R^n$ with $\E^n$ this gives a function which we also write as
$\varphi_\xi : \E^n \to \C$.  The compact group $K$ rotates that
function and we average it over $K$ to obtain
\begin{equation}\label{formula1}
\varphi_\xi^K(x) := \int_K \varphi_{k(\xi)}(x)d\mu_{_K}
= \int_K e^{i b(x,k(\xi))} d\mu_{_K}
\end{equation}
\begin{lemma} \label{formula2}
If $\xi \in \C^n$ then the lift of $\varphi_\xi^K$ from $\E^n$ to $G$
is a $(G,K)$--spherical function.
\end{lemma}
This is standard and can be found in any of \cite{Di}, \cite{H1}, \cite{H2} 
or \cite{W3}.
Curiously, the fact that every $(G,K)$--spherical function is one of the
$\varphi_\xi^K$ seems only to be in the literature for the extreme cases
mentioned above.  We give a general proof below.

In the special case $K = \{1\}$ we have $\varphi_\xi^K = \varphi_\xi$,
resulting in all the quasicharacters on $\R^n$ and thus all the spherical
functions on $\E^n$.

In the trivial case $n = 1$, either $K = \{1\}$ and the matter is described 
above, or $K = \{\pm 1\}$ and $\cD(G,K)$ is the algebra of polynomials in
$\Delta = -\frac{d^2}{dx^2}$.  In the latter case,
the solutions to $\Delta f = \lambda^2f$ on
$\R$ are the linear combinations of $e^{\pm i\lambda x}$, so the
$K$--invariant ones are the multiples of $\cosh(i\lambda x)$.
Thus, if $K = \{\pm 1\}$ we have $\varphi_\xi^K(x) = \cosh(i\sqrt{\xi}\ x)$.

Now suppose $n > 1$.
In the special case $K = O(n)$ or $SO(n)$, the $\varphi_\xi^K$ are
radial functions.  Here $\cD(G,K)$ is the algebra of polynomials in the 
Laplace--Beltrami operator $\Delta = -\sum \frac{\partial^2}{\partial x_i^2}$,
so the joint eigenvalue of $\varphi_\xi^K$ is specified by its
$\Delta$--eigenvalue.  That eigenvalue is $b(\xi,\xi)$, because
$\Delta(\varphi_\xi^K)(x) = 
\Delta_x \int_K e^{i b(x,k(\xi))} d\mu_{_K} =
\int_K \Delta_x e^{i b(x,k(\xi))} d\mu_{_K} =
\int_K b(k(\xi),k(\xi))e^{i b(x,k(\xi))} d\mu_{_K} =
\int_K b(\xi,\xi)e^{i b(x,k(\xi))} d\mu_{_K} =
b(\xi,\xi)\varphi_\xi^K(x)$.  Note that $\varphi_\xi^K$ is a radial
function and the radial part of $-\Delta$ is 
$\frac{d^2}{dr^2} + \frac{n-1}{r}\frac{d}{dr}$.  
Comparing this with the Bessel equation 
$t^2\frac{d^2f}{dt^2} + t \frac{df}{dt} + (t^2 - \nu^2)f = 0$ 
of order $\nu = \frac{n-2}{2}$, we see from \cite[Chap. III]{Wa} 
that the radial eigenfunctions of $\Delta$ for eigenvalue $\lambda^2$ 
are the multiples of
$$
\begin{aligned}
(\lambda r)^{-\nu}J_\nu(\lambda r) 
&= (\lambda r)^{-\nu}\cdot 
 \frac{(\tfrac{1}{2}\lambda r)^\nu}
	{\Gamma(\nu + \tfrac{1}{2})\Gamma(\tfrac{1}{2})}
  \int_0^\pi \cos(\lambda r \cos \theta)\sin^{2\nu}\theta \ d\theta \\
&= \frac{1}{2^\nu\Gamma(\nu + \tfrac{1}{2})\Gamma(\tfrac{1}{2})}
  \int_0^\pi \cos(\lambda r \cos \theta)\sin^{2\nu}\theta \ d\theta 
\text{ where } \nu = \tfrac{n-2}{2}
\end{aligned}
$$
where $J_\nu$ is the Bessel function of order $\nu$.
Note from the last expression that this is even in $\lambda$.
Given $\xi \in \C^n$ we take $\lambda$ to be either of $\pm\sqrt{b(\xi,\xi)}$.
Taking account of $\varphi_\xi^K(0) = 1$ we have
$$
\begin{aligned}
\varphi_\xi^K(x) &= c(n)(\lambda ||x||)^{-(n-2)/2}
       J_{(n-2)/2}(\lambda ||x||) \\
&= \frac{c(n)}{2^{(n-2)/2}
	\Gamma(\tfrac{n-1}{2})\Gamma(\tfrac{1}{2})}
  \int_0^\pi \cos(\lambda ||x|| \cos \theta) \sin^{n-2}\theta \ d\theta \\
&\text{where } \lambda = \sqrt{b(\xi,\xi)} \geqq 0 \text{ and }
	c(n)^{-1}  =  \frac{1}{2^{(n-2)/2}
	\Gamma(\tfrac{n-1}{2})\Gamma(\tfrac{1}{2})}
	  \int_0^\pi \sin^{n-2}\theta \ d\theta.
\end{aligned}
$$
Since $\int_0^\pi \sin^{n-2}\theta \ d\theta =
\sqrt{\pi}\ \Gamma(\frac{n-1}{2})/\Gamma(\frac{n}{2})$ for $n > 1$
we have $c(n) = 
\pi^{-1/2}2^{(n-2)/2} \Gamma(\frac{1}{2}) \Gamma(\frac{n}{2})$, and
\begin{equation}\label{sphon}
\varphi_\xi^K(x) = \frac{\pi^{-1/2} \Gamma(\frac{n}{2})} 
 {\Gamma(\tfrac{n-1}{2})}
 \int_0^\pi \cos(\sqrt{b(\xi,\xi)} ||x||\cos \theta)\sin^{n-2}\theta \ d\theta
  \text{ for } n > 1.
\end{equation}
Of course in general one cannot expect an explicit formula such as 
(\ref{sphon}) for arbitrary cases of $K$, but there is a certain structure,
and now we examine it.

\section{General Spherical Functions on Euclidean Space} \label{sec4}
\setcounter{equation}{0}

In this section, $K$ is an arbitrary closed subgroup of the orthogonal 
group $O(n)$, and $G$ is the semidirect product $\R^n\cdot K$ consisting of
all translations of $\E^n$ with those rotations that are given by elements
of $K$.

As a closed subgroup of the Lie group $O(n)$, $K$ is a compact linear
Lie group.  Thus we can define its complexification $K_{_\C}$  as follows.
The identity component $K_{_\C}^0$ is the analytic subgroup of
$GL(n;\C)$ with Lie algebra $\gk_{_\C}$.  Since
$Ad(K)\gk = \gk$ we have $Ad(\gk_{_\C})
= \gk_{_\C}$, so $K$ normalizes $K_{_\C}^0$.  Thus
$K_{_\C} := KK_{_\C}^0$ is well defined and acts on $\C^n$ as a closed complex
subgroup of the complex orthogonal group $O(n;\C)$.
Note that $K$ is a maximal compact subgroup of $K_{_\C}$.

Recall the categorical quotient $\C^n//K_{_\C}$.  If $\cP(\C^n)^K$ denotes the
algebra of all $/K_{_\C}$--invariant polynomials on $\C^n$, we have the
equivalence relation that $\xi \sim \xi'$ if $p(\xi) = p(\xi')$ for all
$p \in \cP(\C^n)^K$.  Then $\C^n//K_{_\C}$ is the space of equivalence classes,
and it has the structure of affine variety for which $\cP(\C^n)^K$ is the 
algebra of rational functions.  See \cite{Bo} for a good quick development
of this material, \cite{Mu} and \cite{Do} for complete treatments.

\begin{lemma} \label{ft}
Fourier transform gives a $K$--equivariant from the space $\cD(\C^n)$
of constant coefficient differential operators on $\R^n$ onto the space 
$\cP(\C^n)$ of polynomials on $\C^n$.  In particular it gives an 
isomorphism of $\cD(G,K)$ onto $\cP(\C^n)^K$.
\end{lemma}
Without $K$--equivariance, this is a standard basic fact from Fourier 
analysis, and the $K$--equivariance is clear because $b$ is
$K$--invariant.

\begin{theorem} \label{maintheorem1}
Let $\xi, \xi' \in \C^n$.  Then the following conditions are equivalent.

{\rm 1.} The $(G,K)$--spherical functions $\varphi_\xi^K = \varphi_{\xi'}^K$.

{\rm 2.} The orbit closure {\rm c}$\ell K_{_\C}(\xi)$ meets 
{\rm c}$\ell K_{_\C}(\xi')$.

{\rm 3.} If $p$ is a $K_{_\C}$--invariant polynomial on $\C^n$ then 
$p(\xi) = p(\xi')$.

{\rm 4.} The vectors $\xi, \xi' \in \C^n$ have the same image under the 
projection $\pi : \C^n \to \C^n//K_{_\C}$ to the categorical quotient.
\end{theorem}

\noindent {\bf Proof.}
First note that $\xi' \mapsto \varphi_{\xi'}^K(x)$
is constant on every $K$--orbit $\C^n$, for if
$\xi' = k'(\xi'')$ then
$$
\begin{aligned}
\varphi_{\xi'}^K(x)
&= \int_K e^{i b (x,kk'(\xi'')}  d\mu_{_K}(k) 
= \int_K e^{i b (k'^{-1}k^{-1}(x),\xi'')}  d\mu_{_K}(k) \\
&= \int_K e^{i b (k^{-1}(x),\xi'')} d\mu_{_K}(k)
= \int_K e^{i b (x,k(\xi''))} d\mu_{_K}(k) = \varphi_{\xi''}^K(x).
\end{aligned}
$$
Second, note that the function $\varphi_{\xi'}^K(x)$ is holomorphic in  
$\xi'$.  Thus $f(\xi') = \varphi_{\xi'}^K(x)$ is a holomorphic
function on the complex manifold $K_{_\C}(\xi)$ that is constant on
each of the totally real submanifolds $K(\xi'), \ \xi' \in K_{_\C}(\xi)$.
If $J$ is the almost complex structure operator on $K_{_\C}(\xi)$ then
the real tangent spaces satisfy $T_{\xi'}(K_{_\C}(\xi'))
= T_{\xi'}(K(\xi')) + J(T_{\xi'}(K(\xi')))$.  Thus $df = 0$ and $f$ is
constant on every topological component of $K_{_\C}(\xi)$.  But
$K(\xi)$ meet every topological component of $K_{_\mathbb C}(\xi)$, and
$f$ is constant on $K(\xi)$, so now $f$ is constant on $K_{_\C}(\xi)$.
As $f$ is continuous it is constant on {\rm }c$\ell(K_{_\C}(\xi))$
Similarly $f$ is constant on {\rm c}$\ell(K_{_\C}(\xi'))$.  Thus,
if {\rm c}$\ell(K_{_\C}(\xi))$ meets {\rm c}$\ell(K_{_\C}(\xi'))$
then $\varphi_{\xi'}^K = \varphi_\xi^K$.  

Conversely suppose that {\rm c}$\ell(K_{_\C}(\xi))$ does not meet
{\rm c}$\ell(K_{_\C}(\xi'))$.  
Then we have distinct points $z, z' \in  \C^n //K_{_\C}$
such that $\text{\rm c}\ell(K_{_\C}(\xi)) = \pi^{-1}(z)$ and
$\text{\rm c}\ell(K_{_\C}(\xi')) = \pi^{-1}(z')$, so there is a rational
function on $\C^n //K_{_\C}$ with value $0$ at $z$ and $1$ at
$z'$.  That function lifts to a $K_{_\C}$--invariant polynomial $p$
on $\C^n$, by Lemma \ref{ft}.  
Here $p(\xi) = 0$ and $p(\xi') = 1$.  The inverse Fourier transform
of multiplication by $p$ is a $K$--invariant constant coefficient differential
operator $D \in \cD(G/K)$.  By construction $D\varphi_\xi^K = 0$ and
$D\varphi_{\xi'}^K = \varphi_{\xi'}^K$, so $\varphi_\xi^K \ne \varphi_{\xi'}^K$.

We have proved that Assertions 1 and 2 are equivalent.  Equivalence of
Assertions 2, 3 and 4 is standard from invariant theory; see \cite{Bo}
or \cite{Do}.
\hfill $\square$

\begin{remark} \label{rem-equal}
{\rm One could also prove Theorem \ref{maintheorem1} by computing the Fourier
transforms of $\varphi_\xi^K$ and $\varphi_{\xi'}^K$ and using Lemma \ref{ft}
to show equivalence of Assertions 1 and 3, and then using invariant theory
for equivalence of Assertions 2, 3 and 4.}
\end{remark}

Now we know exactly when two of the $(G,K)$--spherical functions $\varphi_\xi^K$
are equal, and we prove that every $(G,K)$--spherical function is one of
them.

\begin{theorem} \label{maintheorem2}
If $\varphi : \E^n \to \C$ is $(G,K)$--spherical, then there
exists $\xi \in \C^n$ such that $\varphi = \varphi_\xi^K$.
\end{theorem}

\noindent {\bf Proof.}  Let $\chi : \cD(G,K) \to \C$ denote the joint 
eigenvalue of $\varphi$ as a joint eigenfunction of $\cD(G,K)$.  Under
Lemma \ref{ft} $\chi$ corresponds to an algebra homomorphism
$\widehat{\chi} : \cP(\C^n)^K \to \C$.  We interpret $\widehat{\chi}$ as
an element of Spec\,$\cP(C^n)^K$.  Thus $\widehat{\chi}$ is evaluation
at some point $[\xi]$ of $\C^n//K_{_\C}$.  Let $\xi \in \pi^{-1}[\xi]$.
Now $\varphi_\xi^K$ is a joint eigenfunction of $\cD(G,K)$ that also
has joint eigenvalue $\chi$.  As both $\varphi$ and $\varphi_\xi^K$ are
$(G,K)$--spherical, they are equal.  \hfill $\square$

We now have all the $(G,K)$--spherical functions parameterized in the form
$\varphi_\xi^K$, and we know exactly when two of them are equal.  In
effect, the affine variety $\C^n//K_{_\C}$ serves as their parameter space.
We make this a little bit more precise.

Let $S(G,K)$ denote the space of all $(G,K)$--spherical functions on $\E^n$.
Recall the spherical transform $f \mapsto \widehat{f}$ from 
$L^1(K \backslash G/K)$ to functions on $S(G,K)$ given by
$\widehat{f}(\varphi) = \int_G f(g)\varphi(g^{-1})d\mu_{_G}$.  Then 
$S(G,K)$ carries the weak topology from the maps 
$\widehat{f} : S(G,K) \to \C$, $f \in L^1(K \backslash G/K)$, or equivalently
from the $\widehat{f} : S(G,K) \to \C$ with $f \in C_c(K \backslash G/K)$.
Further, $S(G,K)$ is a locally compact Hausdorff topological space.

\begin{theorem} \label{ratparam}
The map $\varphi_\xi^K \mapsto \pi(\xi)$ is a homeomorphism of
$S(G,K)$ onto $\C^n//K_{_\C}$.  Thus $S(G,K)$ is parametrized topologically 
by the affine variety $\C^n//K_{_\C}$.
\end{theorem}

\section{Positive Definite Spherical Functions on Euclidean Space} \label{sec5}
\setcounter{equation}{0}

For questions of $(G,K)$--harmonic analysis one must know just which 
$(G,K)$--spherical 
functions $\varphi$ are positive definite, in other words satisfy
$$
\sum \varphi(g_j^{-1}g_i)\overline{c_j}c_i \geqq 0
\text{ whenever } m > 0, \{c_1, \dots , c_m\} \subset \C \text{ and }
\{g_1 , \dots , g_m\} \subset G.
$$
Then $\varphi$ leads to a unitary representation $\pi_\varphi$ on a
Hilbert space $\cH_\varphi$ constructed as follows.  Start with the vector
space of all functions $f : G \to \C$ of finite support.  Give it
the $G$--invariant positive semidefinite hermitian form 
$\langle f, f'\rangle = \sum \varphi(g^{-1}g')f(g)\overline{f'(g')}$.
Dividing out the kernel of the form one has a pre--Hilbert space, 
$\cH_\varphi$ is the completion of that space, and $\pi_\varphi$ is the
natural action of $G$ on $\cH_\varphi$.  See \cite[\S (XXII)9]{Di} or 
\cite[\S 8.6]{W3} 
for the details.  

The representation $\pi_\varphi$ is an irreducible unitary representation 
and its space of $K$--fixed vectors is of dimension $1$.   
Such representations are called $(G,K)$--spherical representations of $G$.
Let $u_\varphi$ be a $K$--fixed unit vector.  Then $\cH_\varphi$ 
determines $\varphi$ (or more precisely its lift to $G$) 
by the formula $\varphi(g) = \langle u, \pi_\varphi(g)u \rangle$.
Again see \cite[\S (XXII)9]{Di} or \cite[Chap. 8]{W3} for details.

If $\varphi = \varphi_\xi^K$ we write $\pi_\xi^K$ for $\pi_\varphi$
and $\cH_{\pi_\xi^K}$ for $\cH_\varphi$\,.  

\begin{proposition} \label{induced-function}
Let $\xi \in \C^n$.  Then $\varphi_\xi^K$ is the induced spherical function
$\text{\rm Ind}_{\R^n}^G(\varphi_\xi)$.  
\end{proposition}

\noindent {\bf Proof.}  Apply formula (\ref{indrep}) to $\varphi_\xi$ 
with $Q = \R^n$.  Since $G$ and $\R^n$ are unimodular, it says that the 
induced spherical function is given by 
$$
\text{\rm Ind}_{\R^n}^G(\varphi_\xi)(xk) =  
	\int_K \varphi_\xi(k_1^{-1}(x)) d\mu_{_K}(k_1) =
	\int_K \varphi_{k_1(\xi)}(x)d\mu_{_K}(k_1) = \varphi_\xi^K(x)
	= \varphi_\xi^K(xk)
$$
for $x \in \R^n$ and $k \in K$.
\hfill $\square$
\medskip

We apply the Mackey little group method to $G$ relative to its normal
subgroup $\R^n$: if $\psi$ is an irreducible unitary representation of
$G$ then it can be constructed (up to unitary equivalence) as follows.
If $\chi$ is a unitary character on $\R^n$
let $K_\chi$ denote its $K$--normalizer, so 
$G_\chi := \R^n\cdot K_\chi$ is the $G$--normalizer of $\chi$.  Write
$\widetilde{\chi}$ for the extension of $\chi$ to $G_\chi$ given by
$\widetilde{\chi}(xk) = \chi(x)$; it is a well defined unitary character 
on $G_\chi$.  If $\gamma$ is an irreducible unitary representation of
$K_\chi$ let $\widetilde{\gamma}$ denote its extension of $G_\chi$ given
by $\widetilde{\gamma}(xk) = \gamma(k)$.  Denote $\psi_{\chi,\gamma} = 
\text{\rm Ind}_{G_\chi}^G(\widetilde{\chi} \otimes \widetilde{\gamma})$.
Then there exist choices of $\chi$ and $\gamma$ such that
$\psi = \psi_{\chi,\gamma}$.

\begin{lemma}\label{kfixed}
In the notation just above, $\psi_{\chi,\gamma}$ has a $K$--fixed 
vector if and only if $\gamma$ is the trivial $1$--dimensional representation
of $K_\chi$.  In that case $\psi_{\chi,\gamma} = 
\text{\rm Ind}_{G_\chi}^G(\widetilde{\chi})$ and the $K$--fixed vector 
is given (up to scalar multiple) by $u(xk) = \chi(x)^{-1}$.
\end{lemma}

\noindent {\bf Proof.} The representation space $\cH_\psi$ of
$\psi = \psi_{\chi,\gamma}$ consists of all $L^2$ functions
$f : G \to \cH_\gamma$ such that 
$f(g'x'k') = \gamma(k')^{-1}\cdot \chi(x')^{-1}f(g')$
for $g' \in G, x' \in \R^n \text{ and } k' \in K_\chi$, and $\psi$ acts 
by $(\psi(g)f)(g') = f(g^{-1}g')$.  

Now suppose that $0 \ne f \in \cH_\psi$
is fixed under $\psi(K)$.  If $k' \in K_\chi$ then $\gamma(k')\cdot f(1)
= f(1)$.  If $f(1) = 0$ then $f(G_\chi) = 0$ and $K$--invariance says
$f = 0$, contrary to assumption.  Thus $f(1) \ne 0$ and irreducibility of
$\gamma$ forces $\gamma$ to be trivial.  

Conversely if $\gamma$ is trivial then $f(xk) = \chi(x)^{-1}$ is a 
nonzero $K$--fixed vector in $\cH_\psi$.  And it is the only one, up to
scalar multiple, because any two $K$--fixed vectors must be proportional.
\hfill $\square$

\begin{lemma} \label{subrep}
In the notation above, $\text{\rm Ind}_{G_\chi}^G(\widetilde{\chi})$
is unitarily equivalent to the subrepresentation of
$\text{\rm Ind}_{\R^n}^G(\widetilde{\chi})$ generated by the 
$K$--fixed unit vector $u(xk) = \chi(x)^{-1}$.
\end{lemma}

\begin{theorem}\label{charposdef}
Let $\varphi$ be a $(G,K)$--spherical function.  Then $\varphi$ is
positive definite if and only if it is of the form $\varphi_\xi^K$
for some $\xi \in \R^n$.  Further, if $\xi, \xi' \in \R^n$ then
$\varphi_\xi^K = \varphi_{\xi'}^K$ if and only if $\xi' \in K(\xi)$.
\end{theorem}

\noindent {\bf Proof.}  Let $\xi \in \R^n$.  
The formula (\ref{formula1}) exhibits
$\varphi_\xi^K$ as a limit of non-negative linear combinations of
positive definite functions on $\R^n$, so it is positive definite.

Now let $\varphi$ be a positive definite $(G,K)$--spherical function.
Let $\pi_\varphi$ be the associated irreducible unitary representation,
and $\cH_\varphi$ the representation space, such that there is a
$K$--fixed unit vector $u_\varphi \in \cH_\varphi$ such that 
$\varphi(g) = \langle  u_\varphi , \pi_\varphi(g)u_\varphi \rangle$
for all $g \in G$.  Following the discussion of the Mackey little group
method, and Lemma \ref{kfixed}, we have a unitary character $\chi$ on
$\R^n$ such that $\pi_\varphi$ is unitarily equivalent to
$\text{\rm Ind}_{G_\chi}^G(\widetilde{\chi})$.  Making the identification,
one $K$--fixed unit vector in $\cH_\varphi$ is given by
$u(xk) = \chi(x)^{-1}$.  We have $\xi \in \R^n$ such that 
$\chi = \varphi_\xi$, so now $u(xk) = \varphi_\xi(x)^{-1}$, and we
compute 
\begin{align}
\varphi(g) &= \langle u, \pi_\varphi(g)u\rangle 
	&&\text{ (construction of $\pi_\varphi$)}\notag \\
&= \langle u, \text{\rm Ind}_{G_{\varphi_\xi}}^G(\widetilde{\varphi_\xi})(g)u
	\rangle &&\text{ (Lemma \ref{kfixed})}\notag \\
&= \langle u, \text{\rm Ind}_{\R^n}^G(\varphi_\xi)(g)u \rangle 
	&&\text{ (construction of $u$)}\notag \\
&= \text{\rm Ind}_{\R^n}^G(\varphi_\xi)(g) 
	&&\text{ (Lemma \ref{subrep})}\notag \\
&= \varphi_\xi^K(g) &&\text{ (Proposition \ref{induced-function})}\notag
\end{align}
That completes the proof of the first assertion.

For the second, suppose that $\xi, \xi' \in \R^n$ with $\varphi_\xi^K
= \varphi_{\xi'}^K$.  Then (up to unitary equivalence)
$\text{\rm Ind}_{G_{\varphi_\xi}}^G(\widetilde{\varphi_\xi}) =
\text{\rm Ind}_{G_{\varphi_{\xi'}}}^G(\widetilde{\varphi_{\xi'}})$.
That gives us direct integral decompositions
$$
\int^\oplus_K \varphi_{k(\xi)} d\mu_{_K}(k)
= \text{\rm Ind}_{G_{\varphi_\xi}}^G(\widetilde{\varphi_\xi})|_{\R^n}
= \text{\rm Ind}_{G_{\varphi_{\xi'}}}^G(\widetilde{\varphi_{\xi'}})|_{\R^n}
= \int^\oplus_K \varphi_{k(\xi')} d\mu_{_K}(k)
$$
All our groups are Type I, so it follows that $K(\xi) = K(\xi')$.
\hfill $\square$
\medskip

We now combine Theorems \ref{ratparam} and \ref{charposdef}.  Let
$P(G,K)$ denote the set of all positive definite $(G,K)$--spherical functions
with the subspace topology from $S(G,K)$.  Then $P(G,K)$ is known to be a
locally compact Hausdorff space, whose topology is the subspace  topology 
from $P(G,K) \subset D^*$ where $D^*$ is the closed unit ball in the
dual space $L^1(K\backslash G/K)^*$.  See, for example, \cite[Prop. 9.2.7]{W3}.

\begin{theorem} \label{ratparamreal}
The map $\varphi_\xi^K \mapsto \pi(\xi)$ is a homeomorphism of $P(G,K)$ onto
$\R^n//K$.  Thus $P(G,K)$ is parameterized topologically by the real form
$\R^n//K$ of the complex affine variety $\C^n//K_{_\C}$.
\end{theorem}

\section{The Transitive Case} \label{sec6}
\setcounter{equation}{0}
Our results are especially interesting
when $K$ is transitive on the spheres about $0$ in $\R^n$,
for then we know the $K_{_\C}$--orbits and the spherical functions
explicitly.
\begin{lemma}\label{quadorbits}
If $K$ is transitive on the spheres about $0$ in $\R^n$, $n > 1$, then
the $K_{_\C}$ is transitive on each of the complex affine quadrics
$Q_c = \{\xi \in \C^n \mid b(\xi,\xi) = c\}, \ 0 \ne c \in \C$.
\end{lemma}

\noindent {\bf Proof.}  If $b(\xi,\xi) = c$ then $K_{_\C}(\xi) \subset Q_c$\,.
Multiplication by a nonzero complex number $t$ gives a $K_{_\C}$--equivariant
holomorphic diffeomorphism of $Q_c$ onto $Q_{t^2c}$, so we need only prove 
the lemma for $0 \ne \xi \in \R^n$.  Then $K(\xi)$ has real dimension
$n-1$, by hypothesis, and is a totally real submanifold of the complex
manifold $K_{_\C}(\xi)$.  It follows that $K_{_\C}(\xi)$ has complex
dimension $n-1$.  Thus $K_{_\C}(\xi)$ is open in $Q_c$.  Now every
$K_{_\C}$--orbit in the connected manifold $Q_c$ is open there, and it
follows that $K_{_\C}$ is transitive on $Q_c$.  \hfill $\square$
\medskip

The situation is more complicated for the cone 
$Q_0 = \{\xi \in \C^n \mid b(\xi,\xi) = 0 \text{ and } \xi \ne 0\}$.
Consider the canonical projection $p : \C^n\setminus \{0\} \to \C\P^{n-1}$
Then $Q:=p(Q_0)$ is the standard nondegenerate projective quadric 
$SO(n)/(SO(n-2) \times SO(2) = SO(n;\C)/P$ where $P$ is the appropriate
parabolic subgroup, and it is known (see \cite{On}, \cite{W2}) that 
$K_{_\C}$ is transitive on $Q$ if and only if
(i) $K_{_\C}^0 = SO(n;\C)$, or (ii) $n = 7$ and $K_{_\C} = G_{_{2,\C}}$,
or (iii) $n = 8$ and $K_{_\C} = Spin(7;\C)$.  In case (i) it is obvious that
$K_{_\C}$ is transitive on $Q_0$\,,  but the argument for all three goes
as follows.  The affine variety $K_{_\C}(\xi)$ is not projective, so
the $K_{_\C}$--equivariant holomorphic map
$p: K_{_\C}(\xi) \to p(K_{_\C}(\xi)) \subset C\P^{n-1}$ has fiber of
dimension $1$.  The image has form $K_{_\C}/P$ where $P$ is a parabolic
subgroup whose reductive component has center $\C^*$, in other words
$P = P' \cdot \C^*$ where $P'$ is the derived group $[P,P]$.  Thus
$K_{_\C}(\xi) \cong K_{_\C}/P'$.  By dimension, now, $K_{_\C}(\xi) = Q_0$\,.

Of course it goes without saying that $K_{_\C}$ is transitive on the
remaining orbit, $\{0\}$.

\begin{lemma} \label{invpoly}
If $K$ is transitive on the spheres about $0$ in $\R^n$, and $p$ is 
a $K_{_\C}$--invariant polynomial on $\C^n$, then $p$ is constant 
on every $Q_c$\,, as well as on $Q_0 \cup \{0\}$.
\end{lemma}

\noindent {\bf Proof.} If $c \ne 0$ then $K_{_\C}$ is transitive on $Q_c$
by Lemma \ref{quadorbits}, so $p$ is constant on $Q_c$\,.  Now let 
$x,y \in Q_0$\,.  Then we have sequences $\{x_m\} \to x$ and
$\{y_m\} \to y$ with $x_m, y_m \in Q_{2^{-m}}$\,.  As $p$ is continuous
now $p(x) = \lim p(x_m) = \lim p(y_m) = p(y)$.  Thus $p$ is constant on
$Q_0$.  As $0$ is in the closure of $Q_0$ also $p$ is constant on
$Q_0 \cup \{0\})$.  \hfill $\square$

\begin{theorem}  \label{catquottrans}
Suppose that $K$ is transitive on the spheres about 
$0$ in $\R^n$.  Then the categorical quotient $C^n//K_{_\C} \cong \C$,
where the isomorphism is given by $Q_c \mapsto c$ for $c \ne 0$ and
$Q_0 \cup \{0\} \mapsto 0$.  
\end{theorem}

\noindent {\bf Proof.}  The projection $\pi : \C^n \to C^n//K_{_\C}$
maps $Q_c$ to a point because $Q_c$ is a closed $K_{_\C}$--orbit, as
noted in Lemma \ref{invpoly}, and maps $Q_0 \cup \{0\}$ to a point
by Lemma \ref{invpoly}.  Thus the map indicated in the statement of 
the theorem is well defined and bijective.  The isomorphism follows.
\hfill $\square$
\medskip

The cases where $K$ is transitive on the spheres about $0$ in $\R^n$
are known explicitly, as follows.
\begin{equation} \label{sphere-trans}
\begin{aligned}
&n > 1 \text{ and } K = SO(n) \text{ or } O(n) \\
&n = 2m \text{ and } K = SU(m) \text{ or } U(m) \\
&n = 4m \text{ and } K = Sp(m) \text{ or } Sp(m)\cdot U(1)  \text{ or }
        Sp(m)\cdot Sp(1) \\
&n = 7 \text{ and } K \text{ is the exceptional group } G_2 \\
&n = 8 \text{ and } K = Spin(7) \\
&n = 16  \text{ and } K = Spin(9)
\end{aligned}
\end{equation}
Now we look at them in some detail.

\begin{lemma} \label{trans-case-invariant-ops}
If $K$ is transitive on the spheres about $0$ in $\R^n$, then
$\cD(G/K) = \C[\Delta]$, algebra of polynomials in the
Laplace--Beltrami operator $\Delta = -\sum \frac{\partial^2}{\partial x_i^2}$.
\end{lemma}

\noindent {\bf Proof.} It is clear that $\C[\Delta] \subset
\cD(G/K)$.  Now let $D \in \cD(G/K)$ be of order $m$.
Then the $m^{th}$ order symbol of $D$ is a polynomial of pure degree
$m$ constant on spheres about $0$ in $\R^n$, in other words a
multiple c$r^m$ with $m$ even and $r^2 = \sum x_i^2$.  Now $D - c(-\Delta)^{m/2}\in \cD(G/K)$ and $D - c(-\Delta)^{m/2}$ has order $< m$.  By
induction on the order, $D - c(-\Delta)^{m/2} \in \C[\Delta]$, so
we have $D \in \C[\Delta]$.
\hfill $\square$
\medskip

Taking advantage of Theorems \ref{ratparam}, \ref{ratparamreal} and 
\ref{catquottrans} we have

\begin{theorem} \label{trans-case-sph-fns}
Let $K$ be transitive on the spheres about $0$ in $\R^n$, $n > 1$.  Then
the spherical functions on $\E^n$ are parametrized by $\C$
and the positive definite spherical functions are the ones with
real non--negative parameter.  Here $\varphi_\xi^K$ has parameter
$b(\xi,\xi)$, which is its $\Delta$--eigenvalue, and
$($given $n$$)$ $\varphi_\xi^K$ is the same for any choice of group
$K$ listed in $(\ref{sphere-trans})$.  Further,
$\varphi_\xi^K(x) = (||\xi||r)^{- (n-2)/2}J_{(n-2)/2}(||\xi||r)
= \frac{\pi^{-1/2} \Gamma(\frac{n}{2})}
 {\Gamma(\tfrac{n-1}{2})}
 \int_0^\pi \cos(\sqrt{b(\xi,\xi)} ||x||\cos \theta)\sin^{n-2}\theta \ d\theta
  \text{ for } n > 1$ 
where $||\xi|| = \sqrt{b(\xi,\xi)}$, $r = \sqrt{b(x,x)}$ and 
$J_\nu$ is the Bessel function of the first  kind of order $\nu$.
\end{theorem}

\begin{remark} \label{quotient}{\rm
Let $M$ be a connected flat $n$--dimensional riemannian symmetric space.
Then $M = \Gamma \backslash \E^n$ for some discrete subgroup $\Gamma$ of
the group $\R^n$ of translations of $\E^n$; see \cite{W1}.  Let $K$ be
a closed subgroup of $O(n)$ that preserves $\Gamma(0)$.  (On the group level 
means that $K$ normalizes $\Gamma$, and if $K$ is connected that forces
$K$ to centralize $\Gamma$.)  Then the action of $K$ descends to $M$, and
$M$ is the symmetric coset space $G'/K$ where $G'$ is the semidirect
product $(\R^n/\Gamma)\cdot K$.  Now $\cD(G',K) = \cD(G,K)$, so any 
$(G',K)$--spherical function on $M$ lifts to a $(G,K)$--spherical function 
on $\E^n$, in other words the $(G',K)$--spherical function on $M$ are the
push--downs of the $(G,K)$--spherical functions $\varphi_\xi^K$ on
$\E^n$ such that $e^{ib(\xi,\gamma)} = 1$ for every $\gamma \in \Gamma$.
So $S(G',K)$ is the categorical quotient $'\C^n//K$ where $'\C^n$ is given
as follows.  Let $U$ denote the complex span $\Gamma \otimes_\Z \R$ of
$\Gamma$ and $U^\perp$ its $b$--orthocomplement; then 
$'\C^n = (U^\perp \oplus U)/\Gamma$.  Note that $'\C^n$ is the annihilator 
of $\Gamma$ in the complexified dual space of $\R^n$.
}
\end{remark}

\noindent Department of Mathematics \hfill\newline
University of California \hfill\newline
Berkeley, CA 94720--3840, USA \hfill\newline
{\tt jawolf@math.berkeley.edu}

\end{document}